\documentclass[12pt]{amsart}

\usepackage{latexsym,rawfonts}
\usepackage{amsfonts,amsmath,amssymb,amsthm}

\usepackage[plainpages=false]{hyperref}
\usepackage{graphicx}
\RequirePackage{color}

\numberwithin{equation}{section}

\newcommand{\beq}{\begin{equation}}
\newcommand{\eeq}{\end{equation}}
\newcommand{\beqs}{\begin{eqnarray*}}
\newcommand{\eeqs}{\end{eqnarray*}}
\newcommand{\beqn}{\begin{eqnarray}}
\newcommand{\eeqn}{\end{eqnarray}}
\newcommand{\beqa}{\begin{array}}
\newcommand{\eeqa}{\end{array}}

\def\p{\partial }

\def\noo{\noindent}

\def\p{\partial}
\def\noo{\noindent}

\def\phi{\varphi}

\begin{document}

 \title []
{\bf  Answer to the questions of \\ Yanyan Li and Luc Nguyen in arXiv:1302.1603}
 
\vskip10pt

\author{Xu-Jia Wang  }

\vskip20pt
 \address{Xu-Jia Wang: Centre for Mathematics and Its Applications,
Australian National University, Canberra, ACT 0200, Australia.}
\email{Xu-Jia.Wang@anu.edu.au. }

 \date{}

\begin{abstract}
In this note we answer the two questions raised by Y.Y Li and L. Nguyen 
in their note [LN2] below.
\end{abstract}
 
\maketitle

\baselineskip=16.2pt
\parskip=5pt

\centerline{\bf ---1---}
 
In the note [LN2], Y.Y Li and L. Nguyen raised two questions.

\begin{itemize}

\item [{\bf Q1.}]  Whether the maximal radial function is super-harmonic.

\item [{\bf Q2.}] A proof of the property $h\to 0$ as $x\to 0$ for bounded $h$, where 
\newline
\centerline{$h(x)=w(x)-2\log|x|$.\hskip50pt\text{$\empty$} }

\end{itemize}

\vskip10pt

\noindent {\bf Answer to Q1}: Given a lower semi-continuous function $v$ in $B_R(x_0)$, 
the {\it maximal radial function} of $v$ is defined by
$$\tilde v(x)=\inf\{v(y):\ y\in\p B_r(x_0),\ r=d(x, x_0)\}, $$
where $B_r(x_0)$ is the geodesic ball of radius $r$ centered at
$x_0$. For any $r\in (0, R)$, there is a point $x_r\in\p B_r(x_0)$ such that
$\tilde v(x_r)=v(x_r)$. 

In page 2445, line -9,  the paper [TW] contains the statement 
``{\it If $v$ is superharmonic, then $\tilde v$ is also superharmonic}.
This statement should be changed to

\noindent 
``{\it If $v$ is superharmonic, then $\tilde v$ is also superharmonic with respect to 
a rotationally symmetric linear operator in $B_r(x_0)$. At any point $x\in B_R(x_0)$,
the coefficients of the operator are equal to those of the Laplacian at $x_r$.
Note that by the exponential map, the Laplacian operator on a manifold 
in local coordinates
is a linear elliptic operator with variable coefficients}.
 
In [TW], we used a $W^{1,p}$ estimate for super-solutions.
This estimate holds for {\it any} linear elliptic equations. 
We would like to thank Y.Y. Li and L. Nguyen for pointing out 
this inaccuracy in our paper.

\vskip10pt

\noindent {\bf Answer to Q2:}  
This question was already answered in my email of November 14, 2012 to Y.Y. Li, 
which was included at the beginning of Section 4 in [W].
 ``{\it with the convergence in $W^{1,p}$, 
if the function $h$ ($h$ is the function in your note) is locally uniformly bounded, 
then the interior gradient estimate or the Harnack inequality 
(for locally bounded solutions) implies the convergence is locally uniform"}. 

I think if one can understand the proof of $|h|\le C$ in page 2456, then  
one should see immediately $h(x)\to 0$ as $x\to 0$,
by repeating the proof in page 2456 and using the interior gradient estimate.
Let me give the details here. 

For any sequence $x_m\to 0$, 
as in [TW] one makes the rescaling: \  
$x\to x/r_m$\ \  (with $r_m=|x_m|$)\  such that \ \ dist$(0, x_m)=1$. 
Denote $A_r=\{x\ |\ 1-r<\text{dist}(0, x)<1+r\}$ the annulus.
We have shown in Lemma 3.4 [TW] that
$$ \text{(i)} \ \ \ \int_{A_{7/8}} |h|\to 0\ \ \text{as} \ \ m\to\infty. $$
From the proof of Theorem 1.3 (page 2456),
$$ \text{(ii)} \ \hskip30pt  \ |h|\le C\ \ \text{in} \ \ A_{3/4}, \hskip20pt\ $$
uniformly in $m$.
By the interior gradient estimate, we have 
$$ \text{(iii)} \ \hskip20pt \ |Dh|\le C\ \ \text{in}\ \ A_{1/2},\hskip20pt\ $$ 
uniformly in $m$.
From (i), (ii), and (iii), we conclude that $h\to 0$ in $A_{1/4}$, uniformly.
Scaling back, we obtain $h(x)\to 0$ as $x\to 0$.

\vskip10pt

Let me pointed out that the main body of the paper [TW]
is to prove (i). 
From (i), one easily obtains (ii). 
The interior gradient estimate (iii) was proved in other papers.

\vskip10pt

\noindent{\bf Remark 1}. When Y.Y. Li asked me Q2 in December 2012,
I thought the answer was already given in [W] and didn't bother to write more.
I just simply said ``{\it there is no need for further correspondence of this mathematics}''.
For Q1, I am sure Y.Y. Li also knew the answer above.

\vskip10pt

\noindent{\bf Remark 2}. 
I didn't know that Y.Y. Li and L. Nguyen posted their note 
[LN2] on arXiv until Wednesday last week. 
I sent the above explanation to them last Friday but have not yet 
received their response for five days. 
So I assume my explanation is clear to them.

\vskip20pt

\newpage

\centerline{\bf ---2---}

\vskip10pt

Based on the questions raised by Y.Y. Li in his emails and in his note [LN2], 
we need to make the following clarifications for the paper [TW].

\begin{itemize}

\item [(1)] (This one is copied from {\bf Answer to Q1} above).\newline
In page 2445, line -9,  the statement 
``{\it If $v$ is superharmonic, then $\tilde v$ is also superharmonic}.
This statement should be changed to 

\vskip5pt

\noindent 
``{\it If $v$ is superharmonic, then $\tilde v$ is also superharmonic with respect to 
a rotationally symmetric linear operator in $B_r(x_0)$. At any point $x\in B_R(x_0)$,
the coefficients of the operator are equal to those of the Laplacian at $x_r$.
Note that by the exponential map, the Laplacian operator on a manifold 
in local coordinates
is a linear elliptic operator with variable coefficients}.

\vskip5pt

\noindent 
Accordingly, Line 1, page 2454, the sentence ``{\it Noticing that $\tilde v_j$ is superharmonic with respect to the conformal Laplace operator (1.16),}'' 
should be changed to 
``{\it Noticing that $\tilde v_j$ is superharmonic with respect to 
a rotationally symmetric linear elliptic operator,}''

\vskip10pt

\item [(2)] (This one is copied from the note [W] below).\newline
After formula (3.29) in page 2455, add

\vskip5pt

\noindent 
``{\it  where $h(x):= w(x)-2\log |x|=o(1)$ is in the sense}
\beqs
 \lim_{r\to 0} r^{-n}\int_{\{r<|x|<2r\}} |h(x)|dx =0,  \text{''}
\eeqs

\vskip10pt

\item [(3)] (This is from {\bf Answer to Q2} above ).
\newline
In page 2456, line 10 after ``{\it This is a contradiction}''
add the new paragraph

\vskip5pt

``{\it This argument also implies that $h(x)\to 0$ as $x\to 0$. 
Indeed, $\forall\ x_m\to 0$, 
make the above rescaling and denote 
$A_r=\{x\ |\ 1-r<\text{dist}(0, x)<1+r\}$.
By Lemma 3.4, we have
$\int_{A_{7/8}} |h|\to 0$ as $m\to\infty$.
The above paragraph tells that $|h|\le C$.
By the interior gradient estimate, 
$|Dh|\le C$  in $A_{1/2}$. Hence  $h\to 0$ in $A_{1/4}$ uniformly as $m\to \infty$.
Scaling back, we obtain $h(x)\to 0$ as $x\to 0$.}''
 
\end{itemize}

\vskip10pt
\noo{\bf Acknowledegment}. 
I would like to thank Y.Y. Li and L. Nguyen 
for giving me the opportunity to make the above
clarifications for the paper [TW]. 
I am sorry that some parts of the paper was not well written 
and have caused difficulties for some readers to understand the proof.

\end{document}